 \documentclass[a4paper,12pt]{article}

 \usepackage{float}
 \usepackage{color}
 \usepackage[pdftex]{graphicx}
 \usepackage{epsfig}       
 \usepackage{subfigure}    
 \usepackage{latexsym}
 \usepackage{graphicx}
 \usepackage{epstopdf}
 \usepackage{amsfonts}
 \usepackage{mathrsfs}
 \usepackage{amsfonts}
 \usepackage{amssymb,amsthm,amsmath,amscd,fancybox,graphicx}
 \usepackage{amssymb,amsmath,euscript}
 \usepackage{amsmath}
 \usepackage{nccmath}
 \usepackage{color}
 \usepackage{fancyhdr}
 \usepackage{multirow}
 \usepackage[figuresright]{rotating}
 \usepackage{empheq}
 \usepackage{cases}
 \usepackage{indentfirst}
 \usepackage[numbers,sort&compress]{natbib}
 
 \usepackage{bm}

 \newcommand{\ds}{\displaystyle}

 \def \beginproof{\par\noindent {\bf Proof.}\ \ }
 \def \endproof{\hskip .5cm $\Box$ \vskip .5cm}

 \newcommand{\M}{\mathcal{M}}
 \newcommand{\I}{\mathcal{I}}
 \newcommand{\D}{\mathcal{D}}
 \newcommand{\V}{\mathcal{V}}

 \DeclareMathOperator{\Tr}{Tr}
 \DeclareMathOperator{\Exp}{Exp}
 \DeclareMathOperator{\intt}{int}

 \theoremstyle{definition}

\begin{document}

 \newtheorem{property}{Property}[section]
 \newtheorem{proposition}{Proposition}[section]
 \newtheorem{append}{Appendix}[section]
 \newtheorem{definition}{Definition}[section]
 \newtheorem{lemma}{Lemma}[section]
 \newtheorem{theorem}{Theorem}[section]
 \newtheorem{example}{Example}[section]
 \newtheorem{corollary}{Corollary}[section]
 \newtheorem{condition}{Condition}
 \newtheorem{remark}{Remark}[section]
 \newtheorem{assumption}{Assumption}[section]
 \newtheorem{algorithm}{Algorithm}[section]
 \newtheorem{problem}{Problem}[section]

 \medskip

 \begin{center}
 {\large \bf  Interval Optimization Problems on Hadamard manifolds}

 \vskip 0.6cm

 Le Tram Nguyen
 \footnote{E-mail: letram07st@gmail.com.} \\
 Faculty of Mathemactics,\ \\
 The University of Da Nang
 - University of Science and Education \\
 and \\
 Department of Mathematics  \\
 National Taiwan Normal University \\
 Taipei 11677, Taiwan

 \vskip 0.6cm

 Yu-Lin Chang \footnote{E-mail: ylchang@math.ntnu.edu.tw.} \\
 Department of Mathematics  \\
 National Taiwan Normal University \\
 Taipei 11677, Taiwan

 \vskip 0.6cm

 Chu-Chin Hu \footnote{E-mail: cchu@ntnu.edu.tw.} \\
 Department of Mathematics  \\
 National Taiwan Normal University \\
 Taipei 11677, Taiwan

 \vskip 0.6cm

 Jein-Shan Chen
 \footnote{Corresponding author. E-mail: jschen@math.ntnu.edu.tw. The research is
 supported by Ministry of Science and Technology, Taiwan.} \\
 Department of Mathematics  \\
 National Taiwan Normal University \\
 Taipei 11677, Taiwan

 \vskip 0.6cm

 May 2, 2022
 \end{center}

 \medskip
 \noindent
 {\bf Abstract.}\
 In this article, we introduce the interval optimization problems (IOPs) on
 Hadamard manifolds as well as study the relationship between them and the
 interval variational inequalities. To achieve the theoretical results, we
 build up some new concepts about $gH$-directional derivative and $gH$-Gâteaux
 differentiability of interval valued functions and their properties on the
 Hadamard manifolds. The obtained results pave a way to further study on
 Riemannian interval optimization problems (RIOPs).

 \medskip


 \medskip

 \section{Motivations}

 This paper studies a new problem set, which is called the interval optimization
 problems on Hadamard manifolds. First, as below, we elaborate the motivation about
 why we focus on this problem. The variational inequalities have been investigated
 since the dawn of the sixties \cite{KS00} and plenty of results are already established,
 see \cite{FP97, C74} and references therein. From \cite{ZZML15}, it is known that
 the interval optimization problems (IOPs) and interval variational inequalities
 (IVIs) possess a close relationship under some assumptions. In addition, the interval
 programming \cite{ALMV21, BP13, GCMD20, C22, S08, W07, W08, ZZML15} is one of the
 approaches to tackle the uncertain optimization problems, in which an interval is
 used to characterize the uncertainty of a variable. Because the variation bounds
 of the uncertain variables can be obtained only through a small amount of uncertainty
 information, the interval programming can easily handle some optimization problems.

 \medskip

 Nowadays, many important concepts and methods of optimization problems have been
 extended from Euclidean space to Riemannian manifolds, particularly to Hadamard
 manifolds \cite{AMS08, CH15, BFO10, BM12, C22, H13, B20}. In general, a manifold
 has no linear structure, nonetheless, it is locally identified with Euclidean space.
 In this setting, the Euclidean metric is replaced by Riemannian metric, which is
 smoothly varying inner product defined on the tangent space of manifold at each
 point, and the line segments is replaced by minimal geodesics. This means that
 the generalization of optimization problems from Euclidean spaces to Riemannian
 manifolds is very important, for example, some nonconvex problems on Euclidean
 space can be viewed as convex problems on Riemannian manifolds. To our best
 knowledge, there is very limited study on Riemannian interval optimization problems
 (RIOPs) in the literature. In \cite{C22}, the authors studied the KKT conditions
 for optimization problems with interval valued objective functions on Hadamard
 manifolds, which is just a routine extension.

 \medskip

 In this paper, we further investigate the interval optimization problems on Hadamard
 manifolds, and characterize the relationship between them and the interval variational
 inequalities. To achieve the theoretical results, we build up some new concepts about
 $gH$-directional derivative and $gH$-Gâteaux differentiability of interval valued
 functions and their properties on the Hadamard manifolds. The analysis differs from
 the one used in traditional variational inequalities and nonlinear programming problems.
 The obtained results pave a way to further study on Riemannian interval optimization
 problems (RIOPs).

 \medskip

 The paper is organized as follows. In Section 2, we formulate the problem set, introduce
 the notations and recall notions of Rienmannian manifolds, tangent space, geodesically
 convex and exponential mapping. We also recall some background materials regarding
 the set of closed, bounded intervals, $gH$-difference and some properties of interval
 valued functions as well as interval valued functions on Hadamard manifolds. In Section
 3, we study the $gH$-continuity, the $gH$-directional derivative and $gH$-Gâteaux
 differentiability of interval valued functions on Hadamard manifolds. Then, we characterize
 the relationship between the $gH$-directional differentiability and geodesically convex
 of Riemannian interval valued functions. In Section 4, we introduce the RIOPs and the
 necessary and sufficient conditions for efficient points of the RIOPs. Besides, we
 define the Riemannian interval variational inequalities problems (RIVIPs) and establish
 the relationship between them and RIOPs. Finally, we draw a conclusion in Section 5.

 \medskip

 \section{Premilinaries}

 In this section, we review some background materials about Riemannian manifolds with
 a special case, the Hadamard manifolds. In particular, we study the intervals, interval
 valued functions, interval valued functions on Hadamard manifolds. We first recall
 some definitions and properties about Riemannian manifolds, which will be used in
 subsequent analysis. These materials can be found in textbooks on Riemannian geometry,
 such as \cite{P16, D92, J05}.

 \medskip

 Let $\M$ be a Riemannian manifold, we denote by $T_x\M$ the tangent space of $\M$
 at $x\in \M$, and the tangent bundle of $\M$ is denoted by $T\M=\cup_{x\in\M}T_x\M$.
 For every $x, y\in\M$, the Riemannian distance $d(x,y)$ on $\M$ is defined by the
 minimal length over the set of all piecewise smooth curves joining $x$ to $y$.
 Let $\nabla$ is the Levi-Civita connection on Riemannian manifold $\M$,
 $\gamma: I\subset \mathbb{R}\longrightarrow\M$ is a smooth curve on $\M$, a vector
 field $X$ is called parallel along $\gamma$ if $ \nabla_{\gamma'}X=0$, where
 $\gamma'=\dfrac{\partial \gamma(t)}{\partial t}$. We say that $\gamma$ is a geodesic
 if $\gamma'$ is parallel along itself, in this cases $\lVert \gamma'\rVert$ is a
 constant. When $\lVert \gamma'\rVert=1, \gamma$ is said to be normalized. A geodesic
 joining $x$ to $y$ in $\M$ is called minimal if its length equals $d(x,y)$.

 \medskip

 For any $x\in\M$, let $V$ be a neighborhood of $0_x\in T_x\M$, the exponential
 mapping $\exp_{x}:V \longrightarrow \M$ is defined by $\exp_x(v)=\gamma(1)$ where
 $\gamma$ is the geodesic  such that $\gamma(0)=x$ and $\gamma'(0)=v$.
 It is known that the derivative of $\exp_x$ at $0_x\in T_x\M$ is the identity map;
 furthermore, by the Inverse Theorem, it is a local diffeomorphism. The inverse map
 of $\exp_x$ is denoted by $\exp_x^{-1}$. A Riemannian manifold is complete if for
 any $x\in\M$, the exponential map $\exp_x$ is defined on $T_x\M$.
 A simply connected, complete Riemannian manifold of nonpositive sectional curvature
 is called a Hadamard manifold. If $\M$ is a Hadamard manifold, for all $x, y\in\M$,
 by the Hopf-Rinow Theorem and Cartan-Hadamard Theorem (see \cite{J05}), $\exp_x$
 is a diffeomorphism and there exists a unique normalized geodesic joining $x$ to
 $y$, which is indeed a minimal geodesic.

 \medskip

 \begin{example} \textbf{Hyperbolic spaces.}
 We equip $\mathbb{R}^{n+1}$ with the Minkowski product defined by
 \[
 \langle x, y\rangle_{1}=-x_0y_0+\sum\limits_{i=1}^{n}x_iy_i,
 \]
 where $x=(x_0, x_1, \cdots, x_n)$, $y=(y_0, y_1, \cdots, y_n)$; and define
 \[
 \mathbb{H}^n:=\{x\in\mathbb{R}^{n+1} \, | \, \langle x, x \rangle_{1}=-1, x_0>0\}.
 \]
 Then, $\langle \cdot ,  \cdot \rangle_{1}$ induces a Riemannian metric $g$ on
 the tangent spaces $T_p\mathbb{H}^{n}\subset \mathbb{R}^{n+1}$, for all
 $p \in \mathbb{H}^{n}$. The section curvature of $(\mathbb{H}^{n}, g)$ is $-1$
 at every point.
 \end{example}

 \medskip

 \begin{example} \textbf{Manifold of symmetric positive definite matrices (SPD).}
 \label{SPD-example}
 The space of $n\times n$ symmetric positive definite matrices  with real entries,
 denoted by $S^{n}_{++}$, is a Hadamard manifold if it is equipped with the below
 Riemannian metric:
 \[
 g_A(x, Y)=\Tr(A^{-1}XA^{-1}Y), \forall A\in S^{n}_{++}, \quad X, Y\in T_{A}S^{n}_{++}.
 \]
 \end{example}

 \medskip
 For more examples, please refer to \cite{B14}. From now on, through the whole paper,
 when we mention $\M$, it means that $\M$ is a Hadamard manifold.

 \begin{definition}[Totally convex set \cite{C94}]
 A subset $\D \subseteq \M$ is said  totally convex if $\D$ contains every geodesic
 $\gamma_{xy}$ of $\M$, whose end points $x, y$ are in $\D$.
 \end{definition}

 \begin{definition}[Geodesically convex set \cite{C94}]
 A subset $\D \subseteq \M$ is said geodesically convex if $\D$ contains the minimal
 geodesic $\gamma_{xy}$ of $\M$, whose end points $x, y$ are in $\D$.
 \end{definition}

 It is easy to see that both total convexity and geodesic convexity are the generalization
 of convexity in Euclidean space. The total convexity is stronger than geodesic
 convexity, but when the geodesic between any two points are unique, they coincide.

 \medskip
 \begin{example}
 Consider $S^{n}_{++}$ as in Example \ref{SPD-example}. Given $a>0$ and let
 \[
 D_a=\{X\in S^{n}_{++} \, | \, \det{X}=a \},
 \]
 then $D_a$ is a nonconvex subset of $S^{n}_{++}$. In fact, from \cite{V18}, the
 minimal geodesic joining $P, Q\in S^{n}_{++}$ is described by
 \[
 \gamma (t)=P^{1/2}(P^{-1/2}QP^{-1/2})^{t}P^{1/2}, \quad \forall t\in[0, 1].
 \]
 If $P, Q\in D_a$, then for all $t\in [0, 1]$ we have
 \begin{eqnarray*}
 \det(\gamma(t))
 &=& \det(P^{1/2}(P^{-1/2}QP^{-1/2})^{t}P^{1/2})  \\
 &=& \det(P)^{1/2}(\det(P)^{-1/2}\det(Q)\det(P)^{-1/2})^{t}\det(P)^{1/2}  \\
 &=& a^{1-t}a^t  \\
 &=& a.
 \end{eqnarray*}
 This means that
 $\gamma(t)\in D_a$, for all $t\in[0, 1]$, that is, $D_a$ is a geodesically convex
 subset of $S^{n}_{++}$.
 \end{example}

 \medskip

 Following the notations used in \cite{DK94}, let $\mathcal{I}(\mathbb{R})$ be the set
 of all closed, bounded interval in $\mathbb{R}$, i.e.,
 \[
 \mathcal{I}(\mathbb{R})=\{[\underline{a}, \overline{a}] \, | \, \underline{a},
 \overline{a}\in\mathbb{R}, \, \underline{a} \leq \overline{a}\}.
 \]
 The Hausdorff metric $d_H$ on $\mathcal{I}(\mathbb{R})$ is defined by
 \[
 d_H(A,B)=\max\{|\underline{a}-\underline{b}|, |\overline{a}-\overline{b}| \}, \quad
 \forall A=[\underline{a}, \overline{a}], \
 B=[\underline{b}, \overline{b}]\in\I(\mathbb{R}).
 \]
 Then, $(\mathcal{I}(\mathbb{R}), d_H)$ is a complete metric space, see \cite{LBD05}.
 The Minkowski sum and scalar multiplications is given respectively by
 \begin{eqnarray*}
 A+B &=& [\underline{a}+\underline{b}, \overline{a}+\overline{b}], \\
 \lambda A &=&
 \begin{cases}
 [\lambda\underline{a}, \lambda\overline{a}]& \text{ if } \lambda \geq 0, \\
 [\lambda\overline{a}, \lambda\underline{a}]& \text{ if } \lambda< 0.
 \end{cases}
 \end{eqnarray*}
 where $A=[\underline{a}, \overline{a}]$, $B=[\underline{b}, \overline{b}]$.
 Note that, $A-A=A+(-1)A \neq 0$.
 A crucial concept in achieving a useful working definition of derivative for
 interval-valued functions is trying to derive a suitable difference between
 two intervals.

 \medskip

 \begin{definition}[$gH$-difference of intervals \cite{S08}]
 Let $A, B\in\mathcal{I}(\mathbb{R})$. The $gH$-difference between $A$ and $B$
 is defined as the interval $C$ such that
 \[
 C=A-_{gH}B \quad \Longleftrightarrow \quad
 \begin{cases}
 A =B+C \\
 \text{or} \\
 B =A-C.
 \end{cases}
 \]
 \end{definition}

 \medskip

 \begin{proposition}{\cite{S08}}
 For any two intervals $A=[\underline{a}, \overline{a}]$, $B=[\underline{b}, \overline{b}]$,
 the $gH$-difference $C=A-_{gH}B$ always exists and
 \[
 C =\left[ \min\{\underline{a}-\underline{b},\overline{a}-\overline{b} \}, \
 \max\{\underline{a}-\underline{b},\overline{a}-\overline{b} \} \right].
 \]
 \end{proposition}

 \medskip

 \begin{proposition}{\cite{LBD05}} \label{property-dH}
 Suppose that $A, B, C\in\I(\mathbb{R})$. Then, the following properties
 hold.
 \begin{description}
 \item[(a)] $d_H(A, B)=0$ if and only if $A=B$.
 \item[(b)] $d_H(\lambda A, \lambda B)=|\lambda|d_H(A, B)$, for all $\lambda\in\mathbb{R}$.
 \item[(c)] $d_H(A+C, B+C)=d_H(A, B)$.
 \item[(d)] $d_H(A+B, C+D) \leq d_H(A, C)+d_H(B, D)$.
 \item[(e)] $d_H(A, B)= d_H(A-_{gH}B, 0)$.
 \item[(f)] $d_H(A-_{gH}B, A-_{gH}C)=d_H(B-_{gH}A, C-_{gH}A)=d_H(B, C)$.
 \end{description}
 \end{proposition}

 Notice that, for all $A\in \I(\mathbb{R})$, we define $||A||:=d_H(A, 0)$, then $||A||$
 is a norm on $\I(\mathbb{R})$ and $d_H(A, B)=||A-_{gH}B||$. There is no natural
 ordering on $\mathcal{I}(\mathbb{R})$, therefore we need to define it.

 \medskip

 \begin{definition}\label{ordering} \cite{W08}
 Let $A=[\underline{a}, \overline{a}]$ and $B=[\underline{b}, \overline{b}]$ be two
 elements of $\mathcal{I}(\mathbb{R})$. We write $A\preceq B$ if
 $\underline{a}\le \underline{b}$ and $\overline{a}\le \overline{b}$. We write
 $A\prec B$ if  $A\preceq B$ and $A\ne B$.
 Equivalently, $A\prec B$ if and only if one of the following cases holds:
 \begin{itemize}
 \item $\underline{a}<\underline{b}$ and $\overline{a}\le \overline{b}$.
 \item $\underline{a}\le\underline{b}$ and $\overline{a}< \overline{b}$.
 \item $\underline{a}<\underline{b}$ and $\overline{a}< \overline{b}$.
 \end{itemize}
 We write, $A\nprec B$ if none of the above three cases hold. If  neither
 $A\prec B$ nor $B\prec A$, we say that none of $A$ and $B$ dominates the other.
 \end{definition}

 \medskip

 \begin{lemma} \label{property-sets}
 For two elements $A, B, C$ and $D$ of $\mathcal{I}(\mathbb{R})$, there hold
 \begin{description}
 \item[(a)] $A\preceq B \ \Longleftrightarrow \ A -_{gH}B\preceq \textbf{0}$.
 \item[(b)] $A\nprec B \ \Longleftrightarrow \ A-_{gH}B\nprec \textbf{0}$.
 \item[(c)] $A\preceq B \ \Longrightarrow \ A-_{gH}C\preceq B-_{gH}C$,
 \item[(d)] $A\preceq B-_{gH}C \ \Longrightarrow \ B\nprec A+C$.
 \item[(e)] $0\preceq (A-_{gH}B)+(C-_{gH}D)\Longrightarrow 0\preceq (A+C)-_{gH}(B+D)$.
 \end{description}
 \end{lemma}
 \beginproof
 (a)  The proofs of part(a)  can be found in \cite{GCMD20}.

 \medskip
 \noindent
 (b) Let $A=[\underline{a}, \overline{a}]$, $B=[\underline{b}, \overline{b}]$, 
 since $A\nprec B$ then $A=B$, or $\underline{a}>\underline{b}$, or 
 $\overline{a}>\overline{b}$. If $\underline{a}>\underline{b}$, or 
 $\overline{a}>\overline{b}$ then $\max\{\underline{a}-\underline{b}, \overline{a}-\overline{b}\}>0$.
 Thus, there holds $\ A-_{gH}B\nprec \textbf{0}$.
 For the other direction, if $\ A-_{gH}B\nprec \textbf{0}$ then $A=B$ or 
 $\max\{\underline{a}-\underline{b}, \overline{a}-\overline{b}\}>0$. This says that 
 $\underline{a}>\underline{b}$, or $\overline{a}>\overline{b}$, which implies $A\nprec B$.

 \medskip
 \noindent
 (c) Let  $A=[\underline{a}, \overline{a}]$, $B=[\underline{b}, \overline{b}]$,
 and $C=[\underline{c}, \overline{c}]$, it is clear that
 \begin{eqnarray*}
 A-_{gH}C &=& \left[ \min\{\underline{a}-\underline{c}, \overline{a}-\overline{c}\},
   \  \max\{\underline{a}-\underline{c}, \overline{a}-\overline{c}\} \right],  \\
 B-_{gH}C &=& \left[ \min\{\underline{b}-\underline{c}, \overline{b}-\overline{c}\},
   \  \max\{\underline{b}-\underline{c}, \overline{b}-\overline{c}\} \right].
 \end{eqnarray*}
 If $A\preceq B$, then $\underline{a} \leq \underline{b}$ and
 $\overline{a} \leq \overline{b}$, which yield
 \[
 \begin{cases}
 \underline{a}-\underline{c} & \leq \underline{b}-\underline{c} \\
 \overline{a}-\overline{c} & \leq \overline{b}-\overline{c}
 \end{cases}
 \quad \Longrightarrow \quad A-_{gH}C \preceq B-_{gH}C.
 \]

 \noindent
 (d) Assume $B\preceq A+C$, by part(c), we know that
 \[
 B-_{gH}C \preceq (A+C)-_{gH}C=A,
 \]
 which indicates $A= B-_{gH}C$. From the definition of $gH$-difference, we have 
 $B=A+C$ or $C=B-A$.
 If $C=B-A$, by the assumption $B\preceq A+C$, there has
 \[
 B\preceq A+B-A \ \Longrightarrow \ A\in\mathbb{R}.
 \]
 Therefore, $B=A+C$. In other words, there holds
 \[
 A \preceq B-_{gH}C \ \Longrightarrow \ B \nprec A+C,
 \]
 which is the desired result.

 \noindent
 (e) Let $A=[\underline{a}, \overline{a}]$, $B=[\underline{b}, \overline{b}]$,
 $C=[\underline{c}, \overline{c}]$ and $D=[\underline{d}, \overline{d}]$.\\
 Case 1: If
 \[
 \begin{cases}
 A-_{gH}B&=[\underline{a}-\underline{b}, \overline{a}-\overline{b}]\\
 C-_{gH}D&=[\underline{c}-\underline{d}, \overline{c}-\overline{d}]
 \end{cases} or \begin{cases}
 A-_{gH}B&=[\overline{a}-\overline{b}, \underline{a}-\underline{b}]\\
 C-_{gH}D&=[\overline{c}-\overline{d}, \underline{c}-\underline{d}]
 \end{cases},
 \]
 then
 \[
 \begin{cases}
 \underline{a}-\underline{b}+\underline{c}-\underline{d}&\ge 0\\
 \overline{a}-\overline{b}+\overline{c}-\overline{d}&\ge 0
 \end{cases}\Longrightarrow \begin{cases}
 (\underline{a}+ \underline{c})-(\underline{b}+\underline{d})&\ge 0\\
 (\overline{a}+ \overline{c})-(\overline{b}+\overline{d})&\ge 0.
 \end{cases}
 \]
 Case 2: If
 \[
 \begin{cases}
 A-_{gH}B&=[\overline{a}-\overline{b}, \underline{a}-\underline{b}]\\
 C-_{gH}D&=[\underline{c}-\underline{d}, \overline{c}-\overline{d}]
 \end{cases}\Longrightarrow \begin{cases}
 \overline{a}-\overline{b}+\underline{c}-\underline{d}&\ge 0\\
 \underline{a}-\underline{b}+\overline{c}-\overline{d}&\ge 0,
 \end{cases}
 \]
 together with
 \[
 \begin{cases}
 \underline{a}-\underline{b}&\ge \overline{a}-\overline{b}\\
 \overline{c}-\overline{d}&\ge \underline{c}-\underline{d}
 \end{cases}
 \]
 we have
 \[
 \begin{cases}
 (\underline{a}+ \underline{c})-(\underline{b}+\underline{d})&\ge 0\\
 (\overline{a}+ \overline{c})-(\overline{b}+\overline{d})&\ge 0.
 \end{cases}
 \]
 Case 3: If
 \[
 \begin{cases}
 A-_{gH}B&=[\underline{a}-\underline{b}, \overline{a}-\overline{b}]\\
 C-_{gH}D&=[\overline{c}-\overline{d}, \underline{c}-\underline{d}]
 \end{cases}\Longrightarrow \begin{cases}
 \overline{a}-\overline{b}+\underline{c}-\underline{d}&\ge 0\\
 \underline{a}-\underline{b}+\overline{c}-\overline{d}&\ge 0,
 \end{cases}
 \]
 together with
 \[
 \begin{cases}
 \overline{a}-\overline{b}&\ge \underline{a}-\underline{b}\\
 \underline{c}-\underline{d}&\ge \overline{c}-\overline{d}
 \end{cases}
 \]
 we have
 \[
 \begin{cases}
 (\underline{a}+ \underline{c})-(\underline{b}+\underline{d})&\ge 0\\
 (\overline{a}+ \overline{c})-(\overline{b}+\overline{d})&\ge 0.
 \end{cases}
 \]
 \endproof

 \begin{remark}
 The inverse of Lemma \ref{property-sets}(c)-(d) are not true. To see this,
 taking $A=[1, 2]$,  $B=[0, 5]$ and $C=[-1,3]$, then
 \[
 A-_{gH}C=[-1, 2],\quad  B-_{gH}C=[1, 2].
 \]
 This means that $A-_{gH}C\preceq B-_{gH}C$, but we do not have $A\preceq B$.
 If taking $A=\textbf{0}$, $B=[0,3]$, $C=[1, 2]$ then
 \[
 A+C=[1, 2], \quad B-_{gH}C=[-1, 1],
 \]
 which says $B\nprec A+C$, but we do not have $A\preceq B-_{gH}C$.
 \end{remark}

 Let $\D\subseteq \M$ be a nonempty set, a mapping
 $f: \D\longrightarrow \mathcal{I}(\mathbb{R})$ is called a Riemannian interval
 valued function (RIVF). We write $f(x)=[\underline{f}(x), \overline{f}(x)]$ where
 $\underline{f}, \overline{f}$ are real valued functions satisfy
 $\underline{f}(x)\le \overline{f}(x)$, for all $x\in\M$. Since $\mathbb{R}^n$
 is a Hadamard manifold, an interval valued function (IVF for short)
 $f:U\subseteq\mathbb{R}^n\longrightarrow \mathcal{I}(\mathbb{R})$ is also a RIVF.

 \medskip

 \begin{definition}\cite{W07}
 Let $U\subseteq\mathbb{R}^{n}$ be a convex set. An IVF
 $f: U\longrightarrow\mathcal{I}(\mathbb{R})$ is said to be convex on $U$ if
 \[
 f(\lambda x_1+(1-\lambda)x_2)\preceq \lambda f(x_1)+(1-\lambda)f(x_2),
 \]
 for all $x_1, x_2 \in U$ and $\lambda\in [0, 1]$.
 \end{definition}

 \medskip

 \begin{definition} \cite{GCMD20}
 Let $U\subseteq\mathbb{R}^n$ be a nonempty set. An IVF
 $f:U\longrightarrow \mathcal{I}(\mathbb{R})$ is said to be
 \textit{monotonically increasing} if for all $x, y\in U$ there has
 \[
 x \leq y \  \Longrightarrow \  f(x)\preceq f(y).
 \]
 The function $f$ is said to be \textit{monotonically decreasing} if for all
 $x, y\in U$ there has
 \[
 x \leq y \  \Longrightarrow \  f(y)\preceq f(x).
 \]
 \end{definition}

 It is clear to see that if an IVF is monotonically increasing (or monotonically
 decreasing), if and only if both the real-valued functions $\underline{f}$ and $\overline{f}$
 are monotonically increasing (or monotonically decreasing).

 \medskip

 \begin{definition}\cite{GCMD20}
 Let $\D\subseteq\M$ be a nonempty set. An RIVF
 $f:\D\longrightarrow \mathcal{I}(\mathbb{R})$ is said to be \textit{bounded below}
 on $\D$ if there exists an interval $A\in\mathcal{I}(\mathbb{R})$ such that
 \[
 A\preceq f(x), \quad \forall x\in \D.
 \]
 The function $f$ is said to be \textit{bounded above} on $\D$ if there exists an
 interval $B\in\mathcal{I}(\mathbb{R})$ such that
 \[
 f(x)\preceq B, \quad \forall x\in \D.
 \]
 The function $f$ is said to be \textit{bounded} if it is both bounded below and
 above.
 \end{definition}
 It is easy to verify that if an RIVF $f$ is bounded below (or bounded above) if and only if
 both the real-valued functions $\underline{f}$ and $\overline{f}$ are bounded
 below (or bounded above).

 \medskip

 \begin{definition}\label{geodesically_convex_set}
 Let $\mathcal{D}\subseteq \M$ be a geodesically convex set and 
 $f:\D\longrightarrow\I(\mathbb{R})$ be a RIVF. $f$ is called geodesically 
 convex on $\D$ if
 \[
 f(\gamma(t))\preceq (1-t) f(x)+tf(y), \quad \forall x, y\in \mathcal{D}
 \ {\rm and} \ \forall t\in [0, 1],
 \]
 where $\gamma:[0, 1]\longrightarrow \M$ is the minimal geodesic joining $x$ and $y$.

 \end{definition}

 \medskip
 \begin{proposition}
 Let $\D$ be a geodesically convex subset of  $\M$ and $f$ be a RIVF on $\D$. 
 Then, $f$ is  geodesically convex on $\D$ if and only if $\underline{f}$ and 
 $\overline{f}$ are geodesically convex on $\D$.
 \end{proposition}
 \beginproof
 This is a direct consequence from Definition \ref{geodesically_convex_set} and 
 Definition \ref{ordering}.
 \endproof

 \medskip

 \begin{example}
 Consider the set
 \[
 \D=\{A\in S^{n}_{++} \, | \, \det(A)>1\}.
 \]
 \end{example}
 \noindent
 For all $X, Y\in \D$, we have the minimal geodesic joining $X, Y$ defined by
 \[
 \gamma(t)=X^{1/2}(X^{-1/2}YX^{-1/2})^{t}X^{1/2}, \quad \forall X, Y\in D
 \ {\rm and} \ \forall t\in [0, 1].
 \]
 For any $t\in [0, 1]$, we also obtain
 \[
 \det(X^{1/2}(X^{-1/2}YX^{-1/2})^{t}X^{1/2})=(\det(X))^{1-t}(\det(Y))^{t}>1,
 \]
 which says that $\D$ is a geodesically convex subset of $S^{n}_{++}$. Moreover,
 on set $\D$, we define a RIVF as below:
 \begin{align}
 f:& \, D\longrightarrow \mathcal{I}(\mathbb{R}) \nonumber \\
   & X\longmapsto \left[0, \, \ln(\det(X)) \right] \nonumber
 \end{align}
 Then, for any $X, Y\in \D$ and $t\in [0, 1]$, we have
 \begin{eqnarray*}
 f(\gamma(t))
 &=& \left[0, \, \ln\det(X^{1/2}(X^{-1/2}YX^{-1/2})^{t}X^{1/2}) \right] \\
 &=& \left[0, \, (1-t)\ln(\det(X))+t\ln(\det(Y)) \right] \\
 &=& (1-t) \left[0, \, \ln(\det(X)) \right] + t \left[0, \, \ln(\det(Y)) \right] \\
 &=& (1-t)f(X) + tf(Y),
 \end{eqnarray*}
 which shows that $f$ is a geodesically convex RIVF on $\D$.

 \medskip
 \begin{proposition}
 The RIVF $f:\D\longrightarrow\I(\mathbb{R})$ is geodesically convex if and only if
 for all $x, y\in\D$ and $\gamma:[0, 1]\longrightarrow \M$ is the minimal geodesic 
 joining $x$ and $y$, the IVF $f\circ\gamma$ is convex on $[0, 1]$.
 \end{proposition}
 \beginproof
 Assume $f$ is geodesically convex, for all $x, y\in \D$ if $\gamma:[0, 1]\longrightarrow \M$ 
 is the minimal geodesic joining $x$ to $y$, then the restriction of $\gamma$ to 
 $[t_1, t_2], t_1, t_2\in [0, 1]$ joins the points $\gamma(t_1)$ to $\gamma(t_2)$. 
 We re-parametrize this restriction
 \[
 \alpha(s)=\gamma(t_1+u(t_2-t_1)), s\in[0, 1].
 \]
 Since $f$ is geodesically convex, for all $s\in [0, 1]$, we have
 \begin{eqnarray*}
 & & f(\alpha(s))\preceq (1-s)f(\alpha(0))+sf(\alpha(1)) \\
 & \Rightarrow & f(\gamma( (1-s)t_1+st_2))\preceq (1-s)f(\gamma(t_1))+sf(\gamma(t_2))  \\
 & \Rightarrow & (f\circ\gamma)( (1-s)t_1+st_2)\preceq (1-s)(f\circ\gamma)(t_1)+s(f\circ\gamma)(t_2),
 \end{eqnarray*}
 which says the IVF $f\circ\gamma$ is convex on $[0, 1]$.
 
 \medskip
 \noindent
 Conversely, for all $x, y\in \D$ and $\gamma:[0, 1]\longrightarrow \M$ is the 
 minimal geodesic joining $x$ and $y$, suppose that
 $f\circ\gamma:[0, 1]\longrightarrow \I(\mathbb{R})$ is a convex IVF. In other words, 
 for all $t_1, t_2\in [0, 1]$, there has
 \[
 (f\circ\gamma)((1-s)t_1+st_2)) \preceq (1-s)(f\circ\gamma)(t_1)+s(f\circ\gamma)(t_2), 
 \ \forall s\in[0, 1].
 \]
 Letting $t_1=0$ and $t_2=1$ gives
 \[
 (f\circ\gamma)(s)\preceq (1-s)(f\circ\gamma)(0)+s(f\circ\gamma)(1), \ \forall s\in[0, 1],
 \]
 or
 \[
 f(\gamma(s))\preceq(1-s)f(x)+sf(y), \ \forall s\in[0, 1].
 \]
 Then, $f$ is a geodesically convex RIVF.
 \endproof

 \medskip
 
 \begin{lemma}
 If $f$ is a geodesically convex RIVF on $\D$ and $A$ is an interval, then the sublevel set
 \[
 \D^{A}=\{x\in\D: f(x)\preceq A\},
 \]
 is a geodesically convex subset of $\D$.
 \end{lemma}
 \beginproof
 For all $x, y\in\D^{A}$, there has $f(x)\preceq A$ and $f(y)\preceq A$. Let 
 $\gamma:[0, 1]\longrightarrow \M$ be the minimal geodesic joining $x$ and $y$. 
 For all $t\in[0, 1]$, by the convexity of $f$, we have
 \[
 f(\gamma(t))\preceq (1-t)f(x)+tf(y)\preceq (1-t)A+tA=A.
 \]
 Thus, $f(\gamma(t))\in \D^{A}$ for all $t\in[0, 1]$, which says that $\D^{A}$ is a 
 geodesically  convex subset of $\D$.
 \endproof

 \medskip

 \section{The $gH$-continuity and $gH$-differentiability of Riemannian interval valued functions}

 In this section, we generalize the $gH$-continuous and $gH$-differentiable
 property of interval valued functions to the settings on the Hadamard manifolds.
 The relationship between  $gH$-differentiability and
 geodesically convex property of the RIVFs is also established.

 \medskip

 \begin{definition}
 Let $f:\M\longrightarrow \mathcal{I}(\mathbb{R})$ be a RIVF,
 $x_0\in\M, A=[\underline{a}, \, \overline{a}]\in \mathcal{I}(\mathbb{R})$. We say
 $\lim_{x\rightarrow x_0}f(x)=A$ if for every $\epsilon>0$, there exists
 $\delta>0$ such that, for all $x\in\M$ and $d(x, x_0)<\delta$, there holds
 $d_H(f(x), A)<\epsilon$.
 \end{definition}

 \medskip
 \begin{lemma}
 Let $f:M\longrightarrow \mathcal{I}(\mathbb{R})$ be a RIVF,
 $A=[\underline{a}, \overline{a}]\in \mathcal{I}(\mathbb{R})$. Then,
 \[
 \lim\limits_{x\rightarrow x_0}f(x)=A
 \quad \Longleftrightarrow \quad
 \begin{cases}
 \lim\limits_{x\rightarrow x_0}\underline{f}(x)=\underline{a}, \\
 \lim\limits_{x\rightarrow x_0}\overline{f}(x)=\overline{a}.
 \end{cases}
 \]
 \end{lemma}
 \beginproof
 If $\lim_{x\rightarrow x_0}f(x)=A$, then for every $\epsilon>0$, there exists
 $\delta>0$ such that, for all $x\in\M$ and $d(x, x_0)<\delta$, there have
 \begin{eqnarray*}
 & & d_H(f(x), A)<\epsilon  \\
 & \Longrightarrow & \max\{|\underline{f}(x)-\underline{a}|, \,
     |\overline{f}(x)-\overline{a}|\} < \epsilon  \\
 & \Longrightarrow &
     \begin{cases}
     |\underline{f}(x)-\underline{a}|<\epsilon \\
     |\overline{f}(x)-\overline{a}|<\epsilon
     \end{cases}.
 \end{eqnarray*}
 Consequently, we have
 \[
 \begin{cases}
 \lim\limits_{x\rightarrow x_0}\underline{f}(x)=\underline{a}, \\
 \lim\limits_{x\rightarrow x_0}\overline{f}(x)=\overline{a}.
 \end{cases}
 \]
 On other hand, if we have
 \[
 \begin{cases}
 \lim\limits_{x\rightarrow x_0}\underline{f}(x)=\underline{a}, \\
 \lim\limits_{x\rightarrow x_0}\overline{f}(x)=\overline{a},
 \end{cases}
 \]
 for every $\epsilon>0$, there exists $\delta>0$ such that, for all $x\in\M$
 and $d(x, x_0)<\delta$, there has
 \[
 \begin{cases}
 |\underline{f}(x)-\underline{a}|<\epsilon \\
 |\overline{f}(x)-\overline{a}|<\epsilon
 \end{cases}
 \quad \Longrightarrow \quad d_H(f(x), A)<\epsilon.
 \]
 which says $\lim\limits_{x\rightarrow x_0}f(x)=A$.
 Thus, the proof is complete.
 \endproof

 \medskip
 \begin{remark}
 From Proposition \ref{property-dH}, we know that $d_H(f(x), A)=d_H(f(x)-_{gH}A, \textbf{0})$,
 which yields
 \[
 \lim\limits_{x\rightarrow x_0}f(x)=A \quad \Longleftrightarrow \quad
 \lim\limits_{x\rightarrow x_0}(f(x)-_{gH}A)=\textbf{0}.
 \]
 \end{remark}

 \medskip

 \begin{definition}[$gH$-continuity] \label{gH-continuity}
 Let $f$ be a RIVF on a nonempty open subset $\mathcal{D}$ of $\M$, $x_0\in \D$.
 The function $f$ is said to be $gH$-continuous at $x_0$ if for all $v\in T_{x_0}\M$
 with $\exp_{x_0}v\in \D$, there has
 \[
 \lim\limits_{||v||\rightarrow 0} \left( f(\exp_{x_0}(v))-_{gH} f(x_0) \right)=\textbf{0}.
 \]
 We call $f$ is $gH$-continuous on $\D$ if $f$ is $gH$-continuous at every $x\in \D$.
 \end{definition}

 \medskip
 \begin{remark}
 We point out couple remarks regarding $gH$-continuity.
 \begin{enumerate}
 \item When $\M\equiv\mathbb{R}^{n}$, $f$ become an IVF and $\exp_{x_0}(v)=x_0+v$. 
       In other words, Definition \ref{gH-continuity} generalizes the concept of 
       the $gH$-continuity of the IVF setting, see \cite{G17}.
 \item By Lemma 3.1 and Remark 3.1, we can see that $f$ is $gH$-continuous if and 
       only if $\underline{f}$ and $\overline{f}$ are continuous.
 \end{enumerate}
 \end{remark}

 \medskip
 
 \begin{theorem}
 Let $\D\subseteq \M$ be a geodesically convex set with nonempty interior and 
 $f: \D\longrightarrow \I(\mathbb{R})$ be a geodesically convex RIVF. Then, $f$ 
 is $gH$-continuous on $\intt\D$.
 \end{theorem}
 \beginproof
 Let $x_0\in\intt \D$ and $B(x_0, r)$ be an open ball center at $x_0$ and of 
 sufficient small radius $r$. Choose $A\in\I(\mathbb{R})$ such that the geodesically 
 convex set $\D^{A}=\{x\in \D: f(x)\preceq A\}$ contains $\overline{B}(x_0, r)$. 
 Let $\gamma: [-1, 1]\longrightarrow\M$ be a minimal geodesic in 
 $\overline{B}(x_0, r)$ such that $\gamma(-1)=x_1, \gamma(0)=x_0, \gamma(1)=x_2$. 
 For convenience, we denote $\gamma(t)=x$ where $t=\frac{d(x_0, x)}{r}\in[0, 1]$. 
 By the convexity of $f$, we have
 \[
 f(\gamma(t))\preceq (1-t)f(x_0)+tf(x_2)\preceq (1-t)f(x_0)+tA,
 \]
 which together with Lemma \ref{property-sets} implies
 \begin{equation} \label{theorem3.11}
 f(x)-_{gH}f(x_0) \preceq t(A-f(x_0)).
 \end{equation}
 The minimal geodesic joining $x_1$ and $x$ is the restriction $\gamma(u), u\in[-1, t]$. 
 Setting $u=-1+s(t+1), s\in[0, 1]$, we obtain the re-parametrization
 \[
 \alpha(s)=y(-1+s(t+1)), \ s\in[0, 1].
 \]
 It is clear to see that
 \[
 \alpha(0)=\gamma(-1)=x_1, \ \alpha\left(\dfrac{1}{t+1}\right)=\gamma(0)=x_0, \
 \alpha(1)=\gamma(t)=x.
 \]
 Due to the convexity of $f$, we have
 \[
 f(\alpha(s))\preceq (1-s)f(x_1)+sf(x)\preceq (1-s)A+sf(x), \ \forall s\in[0, 1].
 \]
 Letting $s=\frac{1}{1+t}$ yields
 \[
 f(x_0)\preceq \dfrac{t}{t+1}A+\dfrac{1}{t+1}f(x),
 \]
 which together with Lemma \ref{property-sets} further implies 
 \begin{equation} \label{theorem3.12}
 f(x_0)-_{gH}f(x)\preceq [(tA+f(x)-_{gH}tf(x_0)]-_{gH}f(x)).
 \end{equation}
 From (\ref{theorem3.11}) and (\ref{theorem3.12}), plugging in $t=\frac{d(x_0, x)}{r}$, we obtain 
 $\lim\limits_{x\rightarrow x_0}f(x)=f(x_0)$. Then, the proof is complete.
 \endproof

 \medskip
 
 \begin{definition} \cite{LMWY11}
 Let $\mathcal{D}\subseteq \M$ be a nonempty open set and consider a function
 $f:\mathcal{D}\longrightarrow\mathbb{R}$. We say that $f$ has directional
 derivative at $x\in \mathcal{D}$ in the direction $v\in T_x\M$ if the limit
 \[
 f'(x,v)=\lim\limits_{t \to 0^{+}}\dfrac{f(\exp_{x}(tv))-f(x)}{t}
 \]
 exists, where $f'(x,v)$ is called the directional derivative of $f$ at $x$
 in the direction  $v\in T_x\M$. If $f$ has directional derivative at $x$ in
 every direction $v\in T_x\M$, we say that $f$ is directional differentiable
 at $x$.
 \end{definition}

 \medskip

 \begin{definition}[$gH$-directional differentiability \cite{C22}]
 Let $f$ be a RIVF on a nonempty open subset $\D$ of $\M$. The function $f$ is
 said to have $gH$-directional derivative at $x\in \D$ in direction $v\in T_x\M$,
 if there exists a closed bounded interval $f'(x,v)$ such that the limits
 \[
 f'(x,v)=\lim\limits_{t \to 0^{+}}\dfrac{1}{t}(f(\exp_x(tv))-_{gH}f(x))
 \]
 exists, where $f'(x,v)$ is called the $gH$-directional derivative of $f$ at $x$
 in the direction of $v$. If $f$ has $gH$-directional derivative at $x$ in every
 direction $v\in T_x\M$, we say that $f$ is $gH$-directional differentiable at $x$.
 \end{definition}

 \medskip

 \begin{lemma}{\cite{C22}}
 Let $\mathcal{D}\subseteq \M$ be a nonempty open set and consider a RIVF
 $f: \mathcal{D}\longrightarrow \mathcal{I}(\mathbb{R})$. Then, $f$ has $gH$-directional
 derivative at $x\in\D$ in the direction $v\in T_x\M$ if and only if $\underline{f}$ and
 $\overline{f}$ have directional derivative at $x$ in the direction $v$. Furthermore,
 we have
 \[
 f'(x,v)=\left[\min\{\underline{f}'(x, v), \overline{f}'(x,v)\}, \,
 \max\{\underline{f}'(x, v), \overline{f}'(x,v)\}\right],
 \]
 where $\underline{f}'(x, v)$ and $\overline{f}'(x,v)$ are the directional derivatives
 of $\underline{f}$ and $\overline{f}$ at $x$ in the direction $v$, respectively.
 \end{lemma}

 \medskip
 
 \begin{theorem}\label{Exitence_gH-directional_derivetive}
 Let $\mathcal{D}\subseteq \M$ be a nonempty open geodesically convex set. If 
 $f:\D\longrightarrow \I(\mathbb{R})$ is a geodesically convex RIVF, then at 
 any $x_0\in\D$, $gH$-directional derivative $f'(x_0, v)$ exists for every 
 direction $v\in T_{x_0}\M$.
 \end{theorem}
 
 \medskip
 \noindent
 To prove Theorem \ref{Exitence_gH-directional_derivetive}, we need two Lemmas.

 \begin{lemma}\label{monotonically_increasing}
 Let $\D\subseteq \M$ be a nonempty geodesically convex set and consider a geodesically 
 convex RIVF $f:\D\longrightarrow \mathcal{I}(\mathbb{R})$. Then, 
 $\forall x_{0}\in\D, v\in T_{x_0}\M$, the function 
 $\phi:\mathbb{R}^{+}\backslash\{0\}\longrightarrow\mathcal{I}(\mathbb{R})$, defined by
 \begin{center}
 $\phi(t)=\dfrac{1}{t}(f(\exp_{x_0}(tv))-_{gH}f(x_0))$,
 \end{center}
 for all $t>0$ such that $\exp_{x_0}(tv)\in\D$, is monotonically increasing.
 \end{lemma}
 \beginproof
 For all $t, s$ such that $0\le t\le s$, by the convexity of $f$, for all $\lambda\in[0, 1]$, we have
 \begin{center}
 $f(\exp_{x_0}(\lambda(sv)))\le (1-\lambda)f(x_0)+\lambda f(\exp_{x_0}(sv))$.
 \end{center}
 Since $\frac{t}{s}\in[0, 1]$, there holds
 \begin{center}
 $f(\exp_{x_0}(tv))\le\dfrac{s-t}{s}f(x_0)+\dfrac{t}{s}f(\exp_{x_0}(sv))$,
 \end{center}
 or
 \begin{align}
 &f(\exp_{x_0}(tv)-_{gH}f(x_0))\nonumber\\
 \le&\left[\dfrac{s-t}{s}f(x_0)+\dfrac{t}{s}f(\exp_{x_0}(sv))\right]-_{gH}f(x_0)\nonumber\\
 =&\left[\min\left\{\dfrac{s-t}{s}\underline{f}(x_0)+\dfrac{t}{s}\underline{f}(\exp_{x_0}(sv))-\underline{f}(x_0), \dfrac{s-t}{s}\overline{f}(x_0)+\dfrac{t}{s}\overline{f}(\exp_{x_0}(sv))-\overline{f}(x_0)\right\},\right.\nonumber\\
 &\left.\max\left\{\dfrac{s-t}{s}\underline{f}(x_0)+\dfrac{t}{s}\underline{f}(\exp_{x_0}(sv))-\underline{f}(x_0), \dfrac{s-t}{s}\overline{f}(x_0)+\dfrac{t}{s}\overline{f}(\exp_{x_0}(sv))-\overline{f}(x_0)\right\}\right]\nonumber\\
 =&\left[\min\left\{\dfrac{t}{s}(\underline{f}(\exp_{x_0}(sv))-\underline{f}(x_0)), \dfrac{t}{s}(\overline{f}(\exp_{x_0}(sv))-\overline{f}(x_0))\right\}\right.\nonumber\\
 &\left.\max\left\{\dfrac{t}{s}(\underline{f}(\exp_{x_0}(sv))-\underline{f}(x_0)), \dfrac{t}{s}(\overline{f}(\exp_{x_0}(sv))-\overline{f}(x_0))\right\}\right]\nonumber\\
 =&\dfrac{t}{s}(f(\exp_{x_0}(sv))-_{gH}f(x_0)).\nonumber
 \end{align}
 Then, the proof is complete.
 \endproof
 
 \medskip
 
 \begin{lemma}\label{bounded_below}
 Let $\D\subseteq\M$ be an open geodesically convex set. 
 If $f:\D\longrightarrow\I(\mathbb{R})$ is a geodesically convex RIVF, then for 
 all $x_0\in \D$ and $v\in T_{x_0}\M$, there exists $t_0\in \mathbb{R}$ such that
 $\phi(t)=\dfrac{1}{t}\left(f(\exp_{x_0}(tv))-_{gH}f(x_0) \right)$ is bounded 
 below for all $t\in(0, t_0]$.
 \end{lemma}
 \beginproof
 For all $v\in T_{x_0}\M$, let $\gamma$ be the geodesic such that $\gamma(0)=x_0$ and 
 $\gamma'(0)=v$. Since $\mathcal{D}\subseteq \M$ be a nonempty open geodesically 
 convex set, there exists $t_1, t_2\in\mathbb{R}$ such that $0\in(t_1, t_2)$ and 
 the restriction of $\gamma$ on $[t_1, t_2]$ is contained in $\D$. Let 
 $\lambda\in (0, t_2]$ and fix the point $\gamma(\lambda)$. The restriction of 
 $\gamma$ to $[t_1, \lambda]$ joins $\gamma(t_1)$ and $\gamma(\lambda)$. We can 
 re-parametrize this restriction
 \[
 \alpha(s)=\gamma(t_1+s(\lambda-t_1)), \ s\in [0, 1].
 \]
 Using the convexity of $f$ gives
 \[
 f(\alpha(s))\preceq (1-s)f(\alpha(0))+sf(\alpha(1)) \ \Longrightarrow \ 
 f(\alpha(s))\preceq (1-s)f(\gamma(t_1))+sf(\gamma(\lambda)) .
 \]
 Plugging in $s=\frac{t_1}{t_1-\lambda}$ leads to
 \[
 f(x_0)\preceq\dfrac{\lambda}{\lambda-t_1}f(\gamma(t_1)) +
 \dfrac{-t_1}{\lambda-t_1}f(\gamma(\lambda)).
 \]
 Then, we have
 \[
 (\lambda-t_1)[\underline{f}(x_0), \overline{f}(x_0)]\preceq 
 \left[ -t_1\underline{f}(\gamma(\lambda))+\lambda\underline{f}(\gamma(t_1)), \ -t_1\overline{f}(\gamma(\lambda))+\lambda\overline{f}(\gamma(t_1)) \right]
 \]
 or
 \[
 \left \{
 \begin{array}{l}
 \dfrac{1}{-t_1}(\underline{f}(x_0)-\underline{f}(\gamma(t_1)))
  \leq  \dfrac{1}{\lambda}(\underline{f}(\gamma(\lambda))-\underline{f}(x_0)),  \\
 \dfrac{1}{-t_1}(\overline{f}(x_0)-\overline{f}(\gamma(t_1)))
  \leq  \dfrac{1}{\lambda}(\overline{f}(\gamma(\lambda))-\overline{f}(x_0)).
 \end{array}
 \right.
 \]
 Thus, the proof is complete.
 \endproof

 \medskip
 \noindent
 As below, we provide the proof of Theorem \ref{Exitence_gH-directional_derivetive}:
 \beginproof
 Let any $x_0\in\D, v\in T_{x_0}\M$. Define an IVF 
 $\phi:\mathbb{R}^{+}\backslash\{0\}\longrightarrow \I(\mathbb{R})$ by
 \[
 \phi(t)=\dfrac{1}{t}(f(\exp_{x_0}(tv))-_{gH}f(x_0)).
 \]
 If $\phi(t)=[\underline{\phi}(t), \overline{\phi}(t)]$, by Lemma \ref{monotonically_increasing} 
 and Lemma \ref{bounded_below}, we have both real-valued functions 
 $\underline{\phi}$ and $\overline{\phi}$ are monotonically increasing and bounded 
 below with $t$ enough small. Therefore, the limits 
 $\lim_{t\rightarrow 0^{+}}\underline{\phi}(t)$ and 
 $\lim_{t\rightarrow 0^{+}}\overline{\phi}(t)$ exist or the limit 
 $\lim_{t\rightarrow 0^+}\phi(t)$ exists.
 Thus, the function $f$ has $gH$-directional derivative at $x_0\in\D$ in the 
 direction $v$.
 \endproof
 
 \medskip
 
 \begin{theorem} \label{geodesical-convexity}
 Let $f:\D\longrightarrow \I(\mathbb{R})$ be a $gH$-directional differentiable RIVF.
 If $f$ is geodesically convex on $\D$, then
 \[
 f'(x, \exp_{x}^{-1}y)\preceq f(y)-_{gH}f(x), \quad \forall x, y\in \D.
 \]
 \end{theorem}
 \beginproof
For all $x, y\in \D$ and $t\in (0, 1]$, by the convexity of $f$, we have
 \[
 f(\gamma(t))\preceq tf(y)+(1-t)f(x),
 \]
 where $\gamma:[0, 1]\longrightarrow \M$ is the minimal geodesic joining $x$ and $y$.
 Applying Lemma \ref{property-sets} yields that
 \begin{eqnarray*}
 & & f(\gamma(t))-_{gH}f(x)  \\
 &\preceq & \left[tf(y)+(1-t)f(x)\right]-_{gH}f(x)  \\
 &=& \left[\min\{t\underline{f}(y)+(1-t)\underline{f}(x)-\underline{f}(x),
     t\overline{f}(y)+(1-t)\overline{f}(x)-\overline{f}(x)\},\right.  \\
 & & \max \left.\{t\underline{f}(y)+(1-t)\underline{f}(x)-\underline{f}(x),
     t\overline{f}(y)+(1-t)\overline{f}(x)-\overline{f}(x)\}\right]  \\
 &=& \left[ \min\{t(\underline{f}(y)-\underline{f}(x)), t(\overline{f}(y)-\overline{f}(x))\},
     \, \max\{t(\underline{f}(y)-\underline{f}(x)),
     t(\overline{f}(y)-\overline{f}(x))\}\right]  \\
 &=& t \, [f(y)-_{gH}f(x)].
 \end{eqnarray*}
 Then, we achieve
 \[
 \dfrac{1}{t}\left[f(\gamma(t))-_{gH}f(x)\right]\preceq f(y)-_{gH}f(x),
 \quad \forall x, y\in\D, \ {\rm and} \ t\in (0, 1].
 \]
 As a result, when $t\longrightarrow 0^{+}$, we obtain
 \begin{center}
 $f'(x, \exp_{x}^{-1}y)\preceq f(y)-_{gH}f(x), \quad \forall x, y\in\D.$
 \end{center}
 Thus, the proof is complete.
 \endproof

 \medskip

 \begin{corollary}
 Let $\D\subseteq\M$ be nonempty open geodesically convex set and suppose that
 the RIVF $f:\D\longrightarrow \I(\mathbb{R})$ is $gH$-directional differentiable
 on $\D$. If $f$ is geodesically convex on $\D$, then
 \[
 f(y)\nprec f'(x, \exp_{x}^{-1}y)+f(x), \quad \forall x, y\in \D.
 \]
 \end{corollary}
 \beginproof
 The result follows immediately from Theorem \ref{geodesical-convexity}
 and Lemma \ref{property-sets}.
 \endproof

 \medskip

 \begin{definition}\cite{GCMD20}
 Let $\V$ be a linear subspace of $\mathbb{R}^{n}$. The IVF
 $F:\V\longrightarrow \I(\mathbb{R})$ is said to be linear if
 \begin{description}
 \item[(a)] $F(\lambda v)=\lambda F(v)$, for all $v\in\V$, $\lambda\in\mathbb{R}$; 
            and
 \item[(b)] for all $ v, w\in \V $, either $F(v)+F(w)=F(v+w)$ or none of $F(v)+F(w)$
            and $F(v+w)$ dominates the other.
 \end{description}
 \end{definition}

 \medskip

 \begin{definition}[$gH$-Gâteaux differentiability]
 Let $f$ be a RIVF on a nonempty open subset $\D$ of $\M$ and $x_0\in \D$. The
 function $f$ is called $gH$-Gâteaux differentiable at $x_0$ if $f$ is $gH$-directional
 differentiable at $x_0$ and $f'(x_0, \cdot): T_{x_0}\M\longrightarrow \I(\mathbb{R})$
 is a $gH$-continuous, linear IVF. The $gH$-Gâteaux derivative of $f$ at $x_0$ is
 defined by
 \[
 f_G(x_0)(\cdot):=f'(x_0, \cdot).
 \]
 The function $f$ is called $gH$-Gâteaux differentiable on $\D$ if $f$ is
 $gH$-Gâteaux differentiable at every $x\in \D$.
 \end{definition}

 \medskip

 \begin{example} \label{example-linear-diff}
 Let $\M:=\mathbb{R}^2$ with the standard metric. Then, $\M$ is a flat Hadamard
 manifold. We consider the RIVF given as below:
 \begin{align}
 f: & \, \M \longrightarrow \I(\mathbb{R}) \nonumber \\
    & (x_1, x_2)\longmapsto
      \begin{cases}
      \dfrac{x_1x_2^2}{x_1^4+x_2^2}[1, 2] & \text{ if } (x_1, x_2)\neq (0, 0), \\
      \textbf{0} & \text { otherwise}.
 \end{cases} \nonumber
 \end{align}
 \end{example}
 \noindent
 For all $v=(v_1, v_2)\in T_{(0,0)}\M\equiv \mathbb{R}^{2}$, we compute
 \begin{eqnarray*}
 f'((0,0), v)
 &=& \lim\limits_{t\rightarrow 0^+}\dfrac{1}{t}\left(f( (0, 0)+tv)-_{gH}f((0,0))\right) \\
 &=& \lim\limits_{t\rightarrow 0^+}\dfrac{1}{t}\dfrac{t^3v_1v_2^2}{t^4v_1^4+t^2v_2^2}[1, 2] \\
 &=& \lim\limits_{t\rightarrow 0^+}\dfrac{v_1v_2^2}{t^2v_1^4+v_2^2}[1, 2] \\
 &=& v_1[1,2].
 \end{eqnarray*}
 On the other hand, for all $h=(h_1, h_2)\in \mathbb{R}^2$, we have
 \begin{eqnarray*}
 & & f'((0, 0), v+h)-_{gH}f'((0, 0), v)  \\
 &=& \left[ \min\{ v_1+h_1-v_1, 2(v_1+h_1)-2v_1\}, \,
     \max\{ v_1+h_1-v_1, 2(v_1+h_1)-2v_1\} \right] \\
 &=& \left[ \min\{h_1, 2h_1\}, \, \max\{h_1, 2h_1\} \right]
 \end{eqnarray*}
 which says
 $\lim\limits_{||h||\rightarrow 0}(f'((0, 0), v+h)-_{gH}f'((0, 0), v)=\textbf{0}$.
 In other words, $f'((0,0), \cdot)$ is a $gH$-continuous IVF. Hence, 
 $f'((0,0), \cdot)$ is a linear, $gH$-continuous IVF or $f$ is $gH$-Gâteaux 
 differentiable at $(0, 0)$ and $f_G( (0, 0))(v)=v_1[1, 2]$ for all 
 $v=(v_1, v_2)\in T_{(0,0)}\M$.

 \medskip

 \begin{example}
 We consider a RIVF defined by
 \begin{align}
 f: & \, S^{n}_{++}\longrightarrow \I(\mathbb{R}) \nonumber \\
    & X\longmapsto
    \begin{cases}
    \left[ \ln(\det (X)), \  \ln(\det( X^2)) \right] & \text{ if } \det (X) \geq 1, \\
    \left[ \ln(\det( X^2)), \ \ln(\det (X)) \right] & \text{ otherwise}.
    \end{cases}\nonumber
 \end{align}
 \end{example}
 \noindent
 For all $ v\in T_IS^{n}_{++}\equiv S^n$, where $S^n$ is the space of  $n\times n$ symmetric
 matrices and $I$ is the $n\times n$ identity matrix, by denoting $Y=\exp_{I}(v)$
 for all $t\in (0, 1]$, we have
 \[
 \ln \left( \det (I^{1/2}(I^{-1/2}YI^{-1/2})^tI^{1/2}) \right) = t\ln(\det( Y)),
 \]
 which implies
 \begin{eqnarray*}
 & & f_G(I)(v)  \\
 &=& \lim\limits_{t\rightarrow 0+}\dfrac{1}{t}[f(\exp_{I}(tv))-_{gH}f(I)]  \\
 &=& \lim\limits_{t\rightarrow 0+}\dfrac{1}{t}\left[\min\{t\underline{f}(Y), t\overline{f}(Y)\},
     \max\{t\underline{f}(Y), t\overline{f}(Y)\}\right],  \\
 &=& \left[ \min\{\ln(\det( Y)), 2\ln(\det( Y))\}, \max \{\ln(\det( Y)), 2\ln(\det (Y))\}\right]  \\
 &=& \begin{cases}
     [\ln(\det( Y)), 2\ln(\det( Y))]& \text{ if } \det( Y) \geq 1 \\
     [2\ln(\det( Y)), \ln(\det( Y))]& \text{ otherwise}.
     \end{cases}
 \end{eqnarray*}
 This concludes that $f$ is $gH$-directional differentiable at $I$.

 \medskip
 \noindent
 On the other hand, for all $v\in S^n, \lambda\in\mathbb{R}$, we know that
 \begin{eqnarray*}
 \exp_{I}(\lambda v)
 &=& I^{1/2}\Exp(I^{-1/2}(\lambda v)I^{-1/2})I^{1/2} \\
 &=& \Exp(\lambda v),
 \end{eqnarray*}
 where $\Exp$ denotes to the matrix exponential. From \cite{H15}, we also have
 \begin{eqnarray*}
 & & \det (\exp_{I}(\lambda v))=\det (\Exp(\lambda v))
     =e^{\Tr (\lambda v)}=\left(e^{\Tr v}\right)^{\lambda}  \\
 &\Longrightarrow & \ln(\det (\exp_{I}(\lambda v)))=\lambda\ln(\det (\exp_{I}( v))).
 \end{eqnarray*}
 To sum up, the function $f_G(x)(.)$ is a linear IVF. Moreover, for
 all $v, h\in S^n$, it follows from \cite{H15} that
 \[
 \exp_{I}(v+h)=\Exp(v+h) =\Exp (v).\Exp( h).
 \]
 Thus, we obtain
 \begin{eqnarray*}
 & & \det (\exp_{I}(v+h))=\det\left(\Exp (v)\right).\det(\Exp (h)) \\
 & \Longrightarrow & \ln(\det (\exp_{I}(v+h)))=\ln(\det (\Exp (v)))+\ln(\det(\Exp (h))). \\
 & \Longrightarrow & \lim\limits_{||h||\rightarrow 0}(f_G(I)(v+h)-_{gH}f_G(I)(v)) \\
 & & =\lim\limits_{||h||\rightarrow 0}[\min\{\ln(\det(\Exp (h))), 2\ln(\det(\Exp (h)))\}, \\
 & & \quad \ \max\{\ln(\det(\Exp (h))), 2\ln(\det(\Exp (h)))\}] \\
 & & =0.
 \end{eqnarray*}
 which says that $f_G(I)(\cdot)$ is a $gH$-continuous IVF. Thus, $f_G(I)(\cdot)$ is
 $gH$-Gâteaux differentiable at $I$.

 \medskip

 \begin{remark}
 We point out that the $gH$-Gâteaux differentiability does not imply the
 $gH$-continuity of RIVF. In fact, in Example \ref{example-linear-diff}, the function
 $f$ is $gH$-Gâteaux differentiable at $(0, 0)$, but
 \[
 \lim\limits_{||h||\rightarrow 0}(f(h_1, h_2)-_{gH}f( (0, 0)))
 =\lim\limits_{||h||\rightarrow 0}\dfrac{h_1h_2^2}{h_1^4+h_2^2}[1, 2],
 \]
 does not exist, which indicates that $f$ is not $gH$-continuous at $(0, 0)$.
 \end{remark}

 \medskip

 \section{Interval optimization problems on Hadamard manifolds}

 This section is devoted to building up some theoretical results on the interval
 optimization problems on Hadamard manifolds. To proceed, we introduce the so-called
 ``efficient point" concept, which is parallel to the role of traditional
 ``minimizer".

 \medskip

 \begin{definition}(Efficient point)
 Let $\D\subseteq\M$ be a nonempty set and
 $f:\D\longrightarrow \I(\mathbb{R})$ be a RIVF. A point $x_0\in\D$ is said to
 be an efficient point of the Riemannian interval optimization problem (RIOP):
 \begin{equation}\label{RIOP}
 \min\limits_{x\in\D} f(x)
 \end{equation}
 if $f(x)\nprec f(x_0)$, for all $x\in \D$.
 \end{definition}

 \medskip
 Sice the objective function $f(x)=[\underline{f}(x), \overline{f}(x)]$ in 
 RIOP (\ref{RIOP}) is an interval-valued function, we can consider two 
 corresponding scalar problems for (\ref{RIOP}) as  follows:
 \begin{equation}\label{LRIOP}
 \min\limits_{x\in\D} \underline{f}(x)
 \end{equation}
 and
 \begin{equation}\label{URIOP}
 \min\limits_{x\in\D} \overline{f}(x)
 \end{equation}
 \medskip
 \begin{proposition}
 Consider problems (\ref{LRIOP}) and (\ref{URIOP}).
 \begin{description}
 \item[(a)] If $x_0\in\D$ is an optimal solution of problems (\ref{LRIOP}) and 
            (\ref{URIOP}) simultaneously , then $x_0$ is an efficient point of 
            the RIOP (\ref{RIOP}).
 \item[(b)] If $x_0\in\D$ is an unique optimal solution of problems (\ref{LRIOP}) 
            or (\ref{URIOP}) , then $x_0$ is an efficient point of the RIOP (\ref{RIOP}).
 \end{description}
 \end{proposition}
 \beginproof
 (a) If $x_0\in\D$ is an optimal solution of problems (\ref{LRIOP}) and (\ref{URIOP}) 
 simultaneously, then
 \[
 \begin{cases}
 \underline{f}(x_0)&\le \underline{f}(x)\\
 \overline{f}(x_0)&\le \overline{f}(x)
 \end{cases}, \forall x\in\D\Rightarrow f(x)\nprec f(x_0), \forall x\in\D,
 \]
 or $x_0$ is an efficient point of RIOP (\ref{RIOP}).
 
 \medskip
 \noindent
 (b) If $x_0\in\D$ is an unique optimal solution of problems (\ref{LRIOP}) or 
 (\ref{URIOP}), then
 \[
 \left[\begin{aligned}
 \underline{f}(x_0)<\underline{f}(x)\nonumber\\
 \overline{f}(x_0)<\overline{f}(x)\nonumber
 \end{aligned}\right., \forall x\in\D\backslash\{x_0\},
 \]
 which says $f(x)\nprec f(x_0)$ for all $x\in \D$, or equivalently $x_0$ is an 
 efficient point of the RIOP (\ref{RIOP}).
 \endproof
 
 \medskip
 
 \begin{proposition}
 Consider the RIOP (\ref{RIOP}) with $f(x)=[\underline{f}(x), \overline{f}(x)]$. 
 Given any $\lambda_1, \lambda_2>0$, if $x_0\in\D$ is an optimal solution of  
 the following problem
 \begin{equation}\label{MP}
 \min\limits_{x\in\D} h(x)=\lambda_1\underline{f}(x)+\lambda_2\overline{f}(x),
 \end{equation}
 then $x_0$ is an efficient point of the RIOP (\ref{RIOP}).
 \end{proposition}
 \beginproof
 Assume that $x_0$ is not an efficient point of RIOP (\ref{RIOP}), then there 
 exists $x'\in\D$ such that
 \[
 f(x')\prec f(x_0) \ \Longrightarrow  \ 
 \lambda_1\underline{f}(x')+\lambda_2\overline{f}(x')
 < \lambda_{1}\underline{f}(x_0)+\lambda_2\overline{f}(x_0).
 \]
  This says that $x_0$ is not an optimal solution of (\ref{MP}), which is a contradiction. 
  Thus, $x_0$ is an efficient point of the RIOP (\ref{RIOP}).
 \endproof

 \medskip
 
 \begin{theorem}[Characterization I of efficient point] \label{characterization-I}
 Let $f: \D\longrightarrow \I(\mathbb{R})$ be a RIVF on a nonempty open subset $\D$
 of $\M$ and $x_0\in \D$ such that $f$ is $gH$-directional differentiable at $x_0$.
 \begin{description}
 \item[(a)] If $x_0$ is an efficient point of  RIOP (\ref{RIOP}), then for all $x\in \D$
            \begin{center}
            $f'(x_0, \exp_{x_0}^{-1}x)\nprec \textbf{0}$ or $f'(x_0, \exp_{x_0}^{-1}x)=[a, 0]$ 
            for some $a<0$.
            \end{center}
 \item[(b)] If $\D$ is geodesically convex, $f$ is geodesically convex  on $\D$ and
            \begin{center}
            $f'(x_0, \exp_{x_0}^{-1}x)\nprec \textbf{0} \quad \forall x\in \D$,
            \end{center}
            then $x_0$ is an efficient point of the RIOP (\ref{RIOP}).
 \end{description}
 \end{theorem}
 \beginproof
 For each $x\in\D$, let $v=\exp_{x_0}^{-1}x$. Since $f$ is $gH$-directional 
 differentiable at $x_0$, then
 \begin{center}
 $f'(x_0, \exp_{x_0}^{-1}x)=f'(x_0, v)=\lim\limits_{t\rightarrow 0^{+}}\dfrac{1}{t}(f(\exp_{x_0}^{-1}(tv))-_{gH}f(x_0))$.
 \end{center}
 (a) If $x_0$ is an efficient point of RIOP (\ref{RIOP}), then
 \begin{align}
 &f(\exp_{x_0}(tv))\nprec f(x_0), \forall t>0  \nonumber \\
 \Rightarrow & f(\exp_{x_0}(tv))-_{gH}f(x_0)\nprec 0, \forall t>0 (\text{ by  Lemma \ref{property-sets}}) \nonumber\\
 \Rightarrow & \dfrac{1}{t}(f(\exp_{x_0}(tv))-_{gH}f(x_0))\nprec 0, \forall t>0  \nonumber\\
 \Rightarrow &\left[\begin{aligned}
 & f'(x_0, v)\nprec \textbf{0} \nonumber\\
 & f'(x_0, v)=[a, 0], \text{ for some } a<0 \nonumber \end{aligned}\right. \nonumber \\
 \Rightarrow & \left[\begin{aligned}
 & f'(x_0, \exp_{x_0}^{-1}x)\nprec \textbf{0} \nonumber \\
 & f'(x_0, \exp_{x_0}^{-1}x)=[a, 0], \text{ for some } a<0 \nonumber 
 \end{aligned}\right. \nonumber
 \end{align}
 
 \medskip
 \noindent
 (b) For all $x\in\D$,  by the convexity of $f$ and applying Proposition 3.1,  we have
 \begin{equation}\label{pt2}
 f'(x_0, \exp_{x_0}^{-1}x)\preceq f(x)-_{gH}f(x_0)\Rightarrow f(x)-_{gH}f(x_0)\nprec \textbf{0}.
 \end{equation}
 On the other hand, by Lemma \ref{property-sets}, there has
 \begin{equation}\label{pt3}
 f(x)\nprec f(x_0)\Leftrightarrow f(x)-_{gH}f(x_0)\nprec \textbf{0} .
 \end{equation}
 From (\ref{pt2}) and (\ref{pt3}), it is clear to see that
 \begin{center}
 $f(x)\nprec f(x_0), \forall x\in\D$.
 \end{center}
 Then, $x_0$ is an efficient point of the RIOP(\ref{RIOP}).
 \endproof

 \medskip

 \begin{example}
 Consider the RIOP $\min\limits_{x\in D} f(x)$ with $f$ and $D$ are defined as
 in Example 2.4. For all $X, Y\in D$ we have
 \[
 f'(X, \exp_{X}^{-1}Y)=
 \begin{cases}
 [0, \, \ln(\det (YX^{-1}))] & \text{ if } \det (Y)\ge\det (X), \\
 [\ln(\det (YX^{-1})), \, 0] & \text{ otherwise}.
 \end{cases}
 \]
 Note that, for all $X\in D$, we can find $Y\in D$ such that $\det Y<\det X$, which
 indicates that this RIOP does not have efficient point.
 \end{example}

 \medskip

 \begin{theorem}[Characterization II of efficient point] \label{characterization-II}
 Let $f: \D\longrightarrow \I(\mathbb{R})$ be a RIVF on a nonempty open subset $\D$
 of $\M$ and $x_0\in \D$ such that $f$ is $gH$-Gâteaux differentiable at $x_0$.
 \begin{description}
 \item[(a)] If $x_0$ is an efficient point of the RIOP (\ref{RIOP}), then
            \[
            0 \in f_{G}(x_0)(\exp_{x_0}^{-1}x), \quad \forall x \in\D.
            \]
 \item[(b)] If $\D$ is geodesically convex, $f$ is a geodesically convex RIVF on $\D$  
            and
            \[
            0 \in [\underline{f}_{G}(x_0)(\exp_{x_0}^{-1}x), \overline{f}_{G}(x_0)(\exp_{x_0}^{-1}x)), 
            \quad \forall x \in\D.
            \]
            where $f_{G}(x_0)(\exp_{x_0}^{-1}x)=[\underline{f}_{G}(x_0)(\exp_{x_0}^{-1}x), \overline{f}_{G}(x_0)(\exp_{x_0}^{-1}x)]$, 
            then $x_0$ is an efficient point of the RIOP (\ref{RIOP}).
 \end{description}
 \end{theorem}
 \beginproof
 For all $x\in\D$, letting $v=\exp_{x_0}^{-1}x$ and due to $f$ being $gH$-Gâteaux
 differentiable at $x_0$, the function $f$ has $gH$-directional derivative at
 $x_0$ in direction $v$ and by Theorem 4.1 we have
 \[
 f_{G}(x_0)(\exp_{x_0}^{-1}x)=f'(x_0, v) \nprec \textbf{0} \text{ or } 
 _{G}(x_0)(\exp_{x_0}^{-1}x)=[a, 0] \text{ for some } a<0.
 \]
 $T_{x_0}\M$ is a linear space then $-v\in T_{x_0}\M$.
 Because $f$ is $gH$-Gâteaux differentiable at $x_0$, hence $f_G(x_0)(\cdot)$ is linear. 
 Then, we obtain
 \[
 f_G(x_0)(-v)=-f_G(x_0)(v),
 \]
 Assume $\underline{f}_G(x_0)(v)>0$, then we have
 \begin{align}
 &-\underline{f}_G(x_0)(v)<0 \nonumber \\
 \Rightarrow & \begin{cases}
 & f_G(x_0)(-v)=[-\overline{f}_G(x_0)(v), -\underline{f}_G(x_0)(v)]\prec\textbf{0} \\
 & f_G(x_0)(-v)\ne [a, 0] \text{ for some } a<0
 \end{cases}, \nonumber
 \end{align}
 which is a contradiction or  $\underline{f}_G(x_0)(v) \leq 0$. Thus, we show that 
 $0\in f_G(x_0)(\exp_{x_0}^{-1}x)$.

 \medskip
 \noindent
 For the remaining part, let $x\in\D$, we have
 \[
 0\in [\underline{f}_{G}(x_0)(\exp_{x_0}^{-1}x), \overline{f}_{G}(x_0)(\exp_{x_0}^{-1}x)) 
 \ \Longrightarrow \  f_G(x_0)(\exp_{x_0}^{-1}x)\nprec \textbf{0}
 \ \Longrightarrow \  f'(x_0, \exp_{x_0}^{-1}x)\nprec \textbf{0},
 \]
 which together with the convexity of $f$ and Theorem \ref{characterization-I}
 proves that  $x_0$ is an efficient point of the RIOP (\ref{RIOP}).
 \endproof

 \medskip

 \begin{example}
 Let $\M=\mathbb{R}_{++}:=\{x\in \mathbb{R} \, | \, x>0\}$ be endowed with  the
 Riemannian metric given by
 \[
 \langle u, v\rangle_x=\dfrac{1}{x^2}uv, \quad \forall u, v\in T_xM\equiv \mathbb{R}.
 \]
 Then, it is known that $\M$ is a Hadamard manifold. For all $x\in\M$, $v\in T_x\M$,
 the geodesic $\gamma:\mathbb{R}\longrightarrow \M$ such that $\gamma(0)=x, \gamma'(0)=v$
 is described by
 \[
 \gamma(t)=\exp_{x}(tv)=xe^{(v/x)t} \quad {\rm and} \quad
 \exp_{x}^{-1}y=x\ln\dfrac{y}{x}, \quad \forall y\in \M.
 \]
 We consider the RIOP $\ds\min\limits_{x\in \M} f(x)$ with
 $f: \M\longrightarrow \I(\mathbb{R})$ is defined by
 \[
 f(x)=\left[x, x+\dfrac{1}{x}\right], \quad \forall x\in \M.
 \]
 For all $x\in \M$, $v\in \mathbb{R}$, we compute
 \begin{eqnarray*}
 f'(x, v)
 &=& \lim\limits_{t\longrightarrow 0^{+}}\dfrac{1}{t}(f(\exp_{x}(tv))-_{gH}f(x))) \\
 &=& \lim\limits_{t\longrightarrow 0^{+}}\dfrac{1}{t}\left[\min\left\{x(e^{(v/x)t}-1),
     x(e^{(v/x)t}-1)+\dfrac{1}{x}(e^{(v/x)t}-1)\right\}\right., \\
 & & \left. \max\left\{x(e^{(v/x)t}-1), x(e^{(v/x)t}-1)+\dfrac{1}{x}(e^{(v/x)t}-1)\right\}\right] \\
 &=& \left[\min\left\{v, v-\dfrac{1}{x^2}v \right\}, \max\left\{v, v-\dfrac{1}{x^2}v \right\}\right] \\
 &=& v\left[1-\dfrac{1}{x^2}, 1\right],
 \end{eqnarray*}
 which says that $f$ is $gH$-directional differentiable on $\M$. We can also easily
 verify that $f'(x, \cdot)$ is $gH$-continuous and linear, and hence $f$ is
 $gH$-Gâteaux differentiable on $\M$.

 \medskip
 \noindent
 On the other hand, by the Cauchy-Schwarz inequality, for all $x>0$, we have
 \begin{center}
 $x+\dfrac{1}{x}\ge 2, $ and $x+\dfrac{1}{x}=2\Leftrightarrow x=1$,
 \end{center}
 then
 \begin{center}
 $\left[x, x+\dfrac{1}{x}\right]\nprec [1, 2], \forall x>0$,
 \end{center}
 or $x=1$ is an efficient point of this RIOP.

 \medskip
 \noindent
 Particularly, at $x_0=1\in \M$,  we have
 \begin{eqnarray*}
 f_G(1)(\exp_{1}^{-1}x)
 &=& \left[\min\left\{\ln x, 0\right\}, \,
     \max\left\{\ln x, 0\right\} \right],  \\
 & = & \begin{cases}
[\ln x, 0]& \text{ if } x<1\\
[0, \ln x]& \text{ if } x\ge 1.
\end{cases}
 \end{eqnarray*}
 \end{example}
 \medskip

 \begin{remark}
 Note that there are similar results in \cite[Theorem 3.2 and Theorem 4.2]{GCMD20},
 which are not correct.
 \begin{enumerate}
 \item From Example 4.1 and Example 4.2, we see that, at $x_0\in \D\subseteq \M$, if 
       there exists $x\in \D$ such that $f'(x_0, \exp_{x_0}^{-1}x)=[a, 0]$ 
       for some $a<0$, we still do not have enough conditions to answer the 
       question: is $x_0$ an efficient point?
 \item Theorem \ref{characterization-I} and Theorem \ref{characterization-II} are 
       the generalization of the Euclidean concepts in \cite[Theorem 3.2 and Theorem 4.2]{GCMD20}. 
       We think their statements are not correct as pointed out as above. Hence, we 
       fix their errors and provide correct versions as in Theorem \ref{characterization-I} 
       and Theorem \ref{characterization-II}.
 \end{enumerate}
 \end{remark}
 
 \medskip
 
 The interval variational inequality problems (IVIPs) was introduced by
 Kinderlehrer and Stampacchia \cite{KS00}. There are some relationships between
 the IVIPs and the IOPs. Let $\D$ be a nonempty subset of $\M$ and
 $T:\D\longrightarrow\mathbb{L}_{gH}(T\M, \I(\mathbb{R}))$ be a mapping such that
 $T(x)\in \mathbb{L}_{gH}(T_x\M,  \I(\mathbb{R}))$, where
 $\mathbb{L}_{gH}(T_x\M, \I(\mathbb{R}))$ denotes the space of $gH$-continuous linear
 mapping from $T_x\M$ to $\I(\mathbb{R})$ and
 $\mathbb{L}_{gH}(T\M, \I(\mathbb{R}))=\bigcup_{x\in\M}\mathbb{L}_{gH}(T_x\M,  \I(\mathbb{R}))$.
 Now, we define the Riemannian interval inequality problems (RIVIPs) as follows:
 \begin{description}
 \item[(a)] The Stampacchia Riemannian interval variational inequality problem (RSIVIP)
            is a problem, which to find $x_0\in\D$ such that
            \[
            T(x_0)(\exp_{x_0}^{-1}y)\nprec\textbf{ 0}, \quad \forall y \in\D.
            \]
 \item[(b)] The Minty Riemannian interval variational inequality problem (RMIVIP) is a
            problem to find $x_0\in\D$ such that
            \[
            T(y)(\exp_{x_0}^{-1}y)\nprec \textbf{0}, \quad \forall y\in\D.
            \]
 \end{description}

 \medskip

 \begin{definition}[Pseudomonotone]
 With a mapping $T$ defined as above, we call $T$ is pseudomonotone if for all
 $x, y\in \D, x\neq y$, there holds
 \[
 T(x)(\exp_{x}^{-1}y)\nprec\textbf{ 0} \quad \Longrightarrow \quad 
 T(y)(\exp_{x}^{-1}y)\nprec \textbf{0}.
 \]
 \end{definition}

 \medskip

 \begin{definition} [Pseudoconvex]
 Let $\D\subseteq \M$ be a nonempty geodesically convex set and
 $f:\M\longrightarrow \I(\mathbb{R})$ be a $gH$-Gâteaux differentiable RIVF.
 Then, $f$ is called pseudoconvex if for all $x, y\in \D$, there holds
 \[
 f_G(x)(\exp_{x}^{-1}y)\nprec \textbf{0} \quad \Longrightarrow \quad f(y)\nprec f(x).
 \]
 \end{definition}

 \medskip

 \begin{proposition}
 Let $\D\subseteq \M$ be a nonempty  set and consider a mapping
 $T: \D\longrightarrow \mathbb{L}_{gH}(T\M, \I(\mathbb{R}))$ such that
 $T(x)\in \mathbb{L}_{gH}(T_{x}\M, I(\mathbb{R}))$. If $T$ is pseudomonotone,
 then every solution of the RSIVIP is a solution of the RMIVIP.
 \end{proposition}
 \beginproof
 Suppose that $x_0$ is a solution of the RSIVIP. Then, we know that
 \[
 T(x_0)(\exp_{x_0}^{-1}y)\nprec \textbf{0}, \quad \forall y\in \D,
 \]
 which together with the pseudomonotonicity of $T$ yields
 \[
 T(y)(\exp_{x_0}^{-1}y)\nprec \textbf{0}, \quad \forall y\in\D.
 \]
 Then, $x_0$ is a solution of the RMIVIP.
 \endproof

 \medskip

 It is observed that if a RIVF $f:\D\longrightarrow \I(\mathbb{R})$ is $gH$-Gâteaux 
 differentiable at $x\in\D$, then $f_{G}(x)\in \mathbb{L}_{gH}(T_x\M,  \I(\mathbb{R}))$. 
 It means there are some relationships between the RIOPs and the RIVIPs.
 
 \medskip
 \begin{theorem}
 Let $\D\subseteq\M$ be a nonempty  set, $x_0\in \D$  and
 $f:\D\longrightarrow\I(\mathbb{R})$ be a $gH$-Gâteaux differentiable RIVF at $x_0$. 
 If $x_0$ is a solution of the RIOP (\ref{RIOP}) and $f_{G}(x_0)(\exp_{x_0}^{-1}y) \neq [a, 0]$
 for all $a<0$, then $x_0$ is a solution of the RSIVIP with $T(x_0)=f_{G}(x_0)$.
 \end{theorem}
 \beginproof
 Since $f$ is $gH$-Gâteaux differentiable on an open set containing $\D$, the function
 $f$ is $gH$-directional differentiable at $x_0$. In light of Theorem
 \ref{characterization-I}, it follows that
 \[
 f_G(x_0)(\exp_{x_0}^{-1}y)\nprec \textbf{0}, \quad \forall y\in \D.
 \]
 Consider $T: \D\longrightarrow \mathbb{L}_{gH}(T\M, \I(\mathbb{R}))$ such that
 $T(x_0)=f_G(x_0)$ for all $x\in \D$, it is clear that
 \[
 T(x_0)(\exp_{x_0}^{-1}y)\nprec \textbf{0}, \quad \forall y\in \D,
 \]
 which says that $x_0$ is a solution of the RSIVIP  with $T(x_0)=f_{G}(x_0)$.
 \endproof

 \medskip

 \begin{theorem}
 Let $\D\subseteq\M$ be a nonempty geodesically convex set, $x_0\in\D$ and
 $f:\D\longrightarrow\I(\mathbb{R})$ be a pseudoconvex, $gH$-Gâteaux differentiable
 RIVF at $x_0$. If $x_0$ is a solution of the RSIVIP with
 $T(x_0)=f_G(x_0)$, then $x_0$ is an efficient point of the RIOP (\ref{RIOP}).
 \end{theorem}
 \beginproof
 Let $x_0$ is a solution of the RSIVIP with $T(x_0)=f_G(x_0)$. Suppose that $x_0$ is
 not an efficient point of RIOP (\ref{RIOP}). Then, exists $y\in\D$ such that
 $f(y)\prec f(x_0)$. From the pseudoconvexity of $f$, we have
 \[
 f_G(x_0)(\exp_{x_0}^{-1}y)\prec\textbf{ 0},
 \]
 which is a contradiction. Thus, $x_0 $ is an efficient point of the RIOP (\ref{RIOP}).
 \endproof

 \medskip

 \section{Conclusions}

 In this paper, we study the Riemannian Interval Optimization problems (RIOPs) on 
 Hadamard manifolds, for which we establish the necessary and sufficient conditions 
 of efficient points. Moreover, we introduce a new concept of $gH$-Gâteaux 
 differentiability of the Riemannian interval valued functions (RIVFs), which is 
 the generalization of $gH$-Gâteaux differentiability of the interval valued functions 
 (IVFs). The Riemannian interval variational inequalities problems (RIVIPs) and their
 relationship, as well as the relationship between the RIVIPs and the RIOPs are also
 investigated in this article. Some examples are presented to illustrate the main
 results.
 \\

 In our opinions, the obtained results are basic bricks towards further investigations
 on the Riemannian interval optimization, in particular, when the Riemannian manifolds
 acts as Hadamard manifolds. For future research, we may either study the theory
 for more general Riemannian manifolds or design suitable algorithms to solve the
 RIOPs.

 \medskip


 \end{document}